\documentclass[12pt]{amsart}
\usepackage{amssymb}

\newcommand{\Z}{{\mathbb Z}}

\newcommand{\R}{{\mathbb R}}
\newcommand{\C}{{\mathbb C}}

\title{New examples of locally algebraically integrable bodies}

\author{V.A.~Vassiliev}
\address{Steklov Mathematical Institute of Russian Academy of Sciences \ \ and \ \ National Research University Higher School of Economics}
 \email{vva@mi-ras.ru}

\begin{document}

\subjclass[2010]{44A99; 4403}

\keywords{Integral geometry, lacuna, algebraic function, Newton's lemma XXVIII}

\begin{abstract}
Any compact body in $\R^N$ with smooth boundary  defines a two-valued function on the space of affine hyperplanes: the volumes of two parts into which these hyperplanes cut the body. This function is never algebraic if $N$ is even and is very rarely algebraic if $N$ is odd: all known examples (found essentially by Archimedes \cite{arch} for $N=3$) of bodies defining algebraic volume functions are these of ellipsoids. We demonstrate a large new series of {\em locally} algebraically integrable bodies with algebraic boundaries in the spaces of arbitrary dimensions, that is, of bodies such that the corresponding volume functions coincide with algebraic ones at least in some open domains of the space of hyperplanes intersecting the body. Further, we discuss (also following Archimedes) the extension of this problem to non-compact bodies.
\end{abstract}

\maketitle

\section{Introduction}

Denote by ${\mathcal P}$ the space of all affine hyperplanes in the Euclidean space $\R^N$.
Any compact subset $W \subset \R^N$ with smooth boundary defines a two-valued function $V_W$ on ${\mathcal P},$ whose values at a hyperplane are the volumes of two parts into which this hyperplane cuts the subset $W$. 
\medskip

\noindent
{\bf Definition} (see \cite{Notices}, \cite{VA}). \rm The compact subset $W\subset \R^N$ is {\em algebraically integrable} if the corresponding volume function $V_W$ coincides everywhere in ${\mathcal P}$ with some pair of values of an algebraic function (i.e., there is a non-trivial polynomial in $N+2$ variables vanishing on any set of numbers $(V, a_1, \dots, a_N, b)$ where $V$ is one of volumes cut from the body $W$ by the hyperplane defined by the equation $a_1x_1 + \dots + a_N x_N =b$). $W$ is {\it locally algebraically integrable} in some open domain of ${\mathcal P}$ if the volume function $V_W$ coincides with an algebraic one in this domain. Any connected domain in ${\mathcal P}$ consisting of hyperplanes transversal to $\partial W$ and such that the local algebraicity property is satisfied in it, is called a {\em lacuna} of $W$, cf. \cite{ABG}, \cite{APLT}. The lacunas consisting of hyperplanes not intersecting the set $W$ are called {\em trivial}. 
\medskip

Archimedes \cite{arch} has shown by an explicit calculation that the balls in $\R^3$ are algebraically integrable; this property can be easily proved also for all ellipsoids in odd-dimensional spaces. On contrary, Newton (\cite{Newton}, Lemma XXVIII) has proved that the convex bounded sets with smooth boundaries in $\R^2$ never are algebraically integrable. V.~Arnold (see \cite{A}, problems 1987-14, 1988-13 and 1990-27) has conjectured that there are no other examples of algebraically integrable bodies than the ellipsoids in odd-dimensional spaces. This conjecture was proved in \cite{Vassiliev-2014} for even-dimensional spaces. Many topological and geometrical obstructions to the integrability in the odd-dimensional case were found in \cite{APLT}, \cite{VA}, however it is unclear whether they are strong enough to prohibit all cases other than the ellipsoids. Moreover, these obstructions imply some hard geometric restrictions on the local integrability property, and motivate the conjecture that for most bodies with semialgebraic boundaries of higher degrees there are no lacunas other than the trivial one. Correspondingly, the problem of indicating all such bodies allowing non-trivial lacunas occurs.

In \cite{VA} a series of bodies with semialgebraic smooth boundaries in any $\R^N$, $N \geq 3$, having non-trivial lacunas was presented: it are some tubular neighbourhoods of standard two-dimensional spheres $S^2 \subset \R^3 \subset \R^N$. In \S \ref{tub} below we expand this series, proving that some tubular neighborhoods of standard even-dimensional spheres $S^{2k} \subset \R^{2k+1} \subset \R^N$ with arbitrary $k< N/2$ also satisfy this local algebraicity condition.
In \S \ref{nob} we discuss (also following Archimedes) the version of the same problem for non-bounded bodies.
\medskip

\noindent
{\bf Remarks.}
1. All calculations below are elementary; the essential part of the work consists in guessing the bodies satisfying the desired property. This guess is based on the geometric obstructions to the local integrability proved in \cite{VA}.

2. The local algebraicity property is affine invariant, therefore all bodies affine equivalent to ones considered in \S 2 below also satisfy it.

3. The condition of the semialgebraicity is essential for the lacuna problem, because in the case of only $C^\infty$-boundaries this problem has many stupid solutions. For instance, if $W \subset \R^n$ is an arbitrary bounded body with smooth boundary, then the corresponding cylindrical body $W \times [-1,1] \subset \R^{n+1}$ with smoothened corners $\partial W \times \{\pm 1\}$ defines a locally algebraic volume function on the space of hyperplanes in $\R^{n+1}$ which are sufficiently close to the hyperplane $\R^n \times \{0\}$. On the other hand, in the case of semialgebraic surfaces the integrability property can be interpreted in the terms of the finiteness of certain orbits of the Picard--Lefschetz monodromy groups related with the surface and simply not defines in the non-algebraic case.

\section{Tubular neighborhoods of spheres}
\label{tub}

Let be $N=n+m$, where $n=2k+1$ is odd and $m$ is an arbitrary natural number. Let $x_1, \dots, x_n$ and $y_1, \dots, y_m$
 be the standard coordinates in the Euclidean spaces $\R^n$ and $\R^m$ respectively.
Consider the $(n-1)$-dimensional sphere in $\R^N=\R^n \times \R^m$ defined by the conditions $x_1^2 + \dots + x_n^2 =1$, $y_1= \dots =y_m=0.$ Given a number $\varepsilon \in (0,1)$, define the body $W \subset \R^N$ as a neighborhood of this sphere, whose boundary $\partial W$ satisfies the equation 
\begin{equation}
(x_1^2+\dots +x_n^2 -1)^2 + \left(y_1^2 + \dots + y_m^2\right)= \varepsilon ^2. \label{tnb}
\end{equation}

This body is invariant under the action of the group ${\mbox O}(n) \times {\mbox O}(m)$ of independent orthogonal transformations of the spaces $\R^n$ and $\R^m$, hence also the volume function $V_W$ is constant on any orbit of the action of this group on the space ${\mathcal P}$ of affine hyperplanes in $\R^{n+m}$. It is easy to see that any such orbit not consisting of hyperplanes parallel to or containing the coordinate subspace $\R^n \equiv \R^n \times \{y\} \subset \R^N $ contains a hyperplane given by the equation
\begin{equation}
x_1=ay_1+b, \qquad a \geq 0, \quad b \geq 0. \label{hh2kp1}
\end{equation}
Namely, the hyperplane defined by the equation 
\begin{equation}
\alpha_1 x_1 + \dots + \alpha_n x_n + \gamma_1 y_1 + \dots + \gamma_m y_m = \beta \ ,
\label{eqhy}
\end{equation}
where not all coefficients $a_j$ are equal to $0,$ can be reduced to the form (\ref{hh2kp1}) with 
\begin{equation} 
a = \frac{\sqrt{\gamma_1^2+ \dots +\gamma_m^2}}{\sqrt{\alpha_1^2+ \dots +\alpha_n^2}}, \quad b= \frac{|\beta|}{\sqrt{\alpha_1^2+ \dots + \alpha_n^2}}.
\label{nof}
\end{equation}
The value $b$ here is the distance in the coordinate subspace $\R^n $ between the origin and the intersection of our hyperplane (\ref{eqhy}) with this subspace; $a$ is the cotangent of the angle between this hyperplane and $\R^n$.

\medskip

\noindent
{\bf Notation}. Let $v_p$ be the volume of the unit ball in the Euclidean space $\R^p$ (so that the $(p-1)$-dimensional area of the unit sphere is equal to $p v_p$). Denote by $L(a,b)$ the hyperplane defined by (\ref{hh2kp1}), and by $V(a,b)$ the bigger value of the volume function $V_W$ at this hyperplane $L(a,b)$ (that is, the volume of the part of $W \cap (\R^N \setminus L(a,b))$ containing the origin in $\R^N$). Let $C(\varepsilon)$ be the volume of the entire body $W$.
\medskip

The body $W$ is symmetric, therefore both values of the function $V_W$ at any hyperplane containing the origin are equal to $C(\varepsilon)/2$; in particular $V(a,0) \equiv C(\varepsilon)/2$.
\medskip

\noindent
{\bf Theorem.}
{\it For any $\varepsilon \in (0,1)$, the partial derivative $\frac{\partial V(a,b)}{\partial b}$ is equal to a polynomial of degree $k-1$ in the variables $a^2$ and $b^2$ in the area of $\R^2_+$ where $a$ and $b$ are sufficiently small with respect to $1-\varepsilon$.}
\medskip

\noindent
{\bf Corollary.}
{\it For any $\varepsilon \in( 0, 1)$, the volume function $V_W$ defined by the body $W$ with the boundary $($\ref{tnb}$)$ in a neighborhood of the set of hyperplanes containing the coordinate subspace $\R^m$ is equal 
to a two-valued function of the parameters $\alpha_i, \gamma_i$ and $\beta$ of these hyperplanes $($see $($\ref{eqhy}$))$ having the form
 $$C(\varepsilon)/2 \pm P\left(a(\alpha, \gamma), b(\alpha,\beta)\right) \ ,$$ where the functions
$a(\alpha,\gamma)$ and $b(\alpha,\beta)$ are given in $($\ref{nof}$)$, and 
$P$ is a polynomial of degree $2k-1$ containing only monomials of even degrees in $a$ and odd degrees in $b$. 

In particular, the function $V_W$ coincides in this neighbourhood with an algebraic function of the variables $\alpha_i, \gamma_j$ and $\beta$.}
\medskip
 
\noindent
{\it Proof of Theorem.} The partial derivative of the volume function $V(a,b)$ with respect to $b$ is obviously equal to 
the $(N-1)$-dimensional volume of the orthogonal projection of the hypersurface $L(a,b) \cap W$ 
to the coordinate subspace $\{x_1=0\}$ in $\R^N$. Let us calculate this volume. 

For any $y \in \R^m$ with $|y| \leq \varepsilon $, consider the corresponding affine $n$-dimensional plane $\R^n_y \equiv \R^n \times \{y\} \subset \R^n \times \R^m$. The intersection set of this plane and the body $W$ is the spherical layer defined by the condition $$(x_1^2+ \dots +x_n^2 -1)^2 \leq \varepsilon^2-|y|^2,$$ i.e. the set of points in $\R_y^n$ whose distances from the origin $(0, y)$ of this plane belong to the segment
$\left[\sqrt{1-\sqrt{\varepsilon^2-|y|^2}} , \sqrt{ 1+\sqrt{\varepsilon^2-|y|^2}}\right]$. If the constants $a$ and $b$ in (\ref{hh2kp1}) are sufficiently small, then the intersection of any such layer with the affine plane (\ref{hh2kp1}) is a spherical layer in some affine subspace of dimension $n-1\equiv 2k$ in $\R^n_y$: the latter layer is the difference of two $(n-1)$-dimensional balls with the radii equal to
$$\sqrt{1-(ay_1+b)^2+\sqrt{\varepsilon^2-|y|^2} } \ \mbox{ and } \ \sqrt{ 1-(ay_1+b)^2-\sqrt{\varepsilon^2-|y|^2} }.$$
The $(n-1)$-dimensional Euclidean volume of this layer 
is equal to 
\begin{equation}
v_{2k} \left( \left(1-(ay_1+b)^2+\sqrt{\varepsilon^2-|y|^2}\right)^k  - \left( 1-(ay_1+b)^2-\sqrt{\varepsilon^2-|y|^2}\right)^k \right) \label{subball}
\end{equation}

We need to integrate these volumes over the ball in $\R^m$ consisting of all points $y$ with $|y|\leq \varepsilon$. Let us fiber this ball into the family of $(m-1)$-dimensional balls, on any of which the coordinate $y_1$ is constant. Then we integrate over $y_1 \in [-\varepsilon, \varepsilon]$ the integrals of the values (\ref{subball}) over the corresponding $(m-1)$-dimensional balls of this family, i.e. in the case $m>1$ we take the integral 

\begin{equation} \int_{-\varepsilon}^\varepsilon   \int_{0}^{\sqrt{\varepsilon^2 - y_1^2}} (m-1)v_{m-1} \rho^{m-2} v_{2k} \left( \left(1-(ay_1+b)^2+\sqrt{\varepsilon^2-(y_1^2 + \rho^2)}\right)^k  - \right. $$
$$\left. - \left( 1-(ay_1+b)^2-\sqrt{\varepsilon^2-(y_1^2+\rho^2)}\right)^k \right) d \rho \ dy_1.
\label{desi}
\end{equation}

The substitution $y_1 =\varepsilon \sin \varphi,$ $ \rho = \varepsilon \cos \varphi \sin \vartheta$, turns this integral to 
$$(m-1)v_{m-1} v_{2k}\int_{-\pi/2}^{\pi/2} \int_0^{\pi/2} (\varepsilon \cos \varphi \sin \vartheta)^{m-2}
 \left( \left(1-(a\varepsilon \sin \varphi+b)^2+ \varepsilon \cos \varphi \cos \vartheta \right)^k - \right. $$
$$ \left. -
 \left( 1-(a\varepsilon \sin \varphi +b)^2- \varepsilon \cos \varphi \cos \vartheta \right)^k \right)
 d (\varepsilon \cos\varphi \sin \vartheta) \ d(\varepsilon \sin \varphi) \equiv$$ 
$$\equiv (m-1)v_{m-1} v_{2k} \varepsilon^{m} \int_{-\pi/2}^{\pi/2}\cos^{m} \varphi \int_0^{\pi/2}   
 \left( \left(1-(a\varepsilon \sin \varphi+b)^2+ \varepsilon \cos \varphi \cos \vartheta \right)^k - \right. $$
$$ \left. -
 \left( 1-(a\varepsilon \sin \varphi +b)^2- \varepsilon \cos \varphi \cos \vartheta \right)^k \right) \sin^{m-2}\vartheta \cos \vartheta \ d \vartheta \ d \varphi.$$

The function under the interior integral in the last expression is a polynomial in $\sin \vartheta$ and $\cos \vartheta$, whose coefficients are some polynomials in the variables $\varepsilon, \cos \varphi, b,$ and $a \sin \varphi$. Therefore the values of these integrals corresponding to different values of $\varphi$ also form a polynomial in these variables. Any monomial of the latter polynomial containing an odd power of $a$ contains also the same odd power of $\sin \varphi$ and hence vanishes in the final integration over $\varphi \in \left[-\frac{\pi}{2}, \frac{\pi}{2}\right]$. The result of this integration is a polynomial in $a$ and $b$ (for any fixed $\varepsilon$) containing only even powers of $a$. But the sum of powers of $a$ and $b$ in any monomial of (\ref{subball}) also is even, hence the same is true for the final polynomial, which therefore contains only even powers of $b$, too. 

In the remaining case of $m=1$, the function (\ref{desi}) is equal to the integral
\begin{equation}
v_{2k} \int_{-\varepsilon}^{\varepsilon} \left(1-(ay+b)^2 + \sqrt{\varepsilon^2-y^2}\right)^k - \left(1-(ay+b)^2 - \sqrt{\varepsilon^2-y^2}\right)^k dy =
\end{equation}
$$
= v_{2k} \varepsilon \int_{-\pi/2}^{\pi/2} \left(1-(a \varepsilon \sin \varphi +b)^2 + \varepsilon \cos \varphi\right)^k - \left(1-(a \varepsilon \sin \varphi +b)^2 - \varepsilon \cos \varphi\right)^k \cos \varphi d \varphi =
$$
$$
= 2 v_{2k} \sum_{j=1}^{]k/2[} \binom{k}{2j-1}\  \varepsilon^{2j} \int_{-\pi/2}^{\pi/2} \left(1-(a\varepsilon \sin \varphi +b)^2\right)^{k-2j+1} \cos^{2j} \varphi \  \ d \varphi ,
$$
also satisfying the assertion of our theorem by the same reasons as above.
\hfill $\Box$

\section{The algebraic integrability of non-bounded domains (after Archimedes)}
\label{nob}

Let $W$ be a non-compact subset with semialgebraic boundary in $\R^N$. It may happen that the volume of one of two parts, into which a hyperplane in $\R^N$ cuts $W$, is finite (and hence the same is true for all neighboring hyperplanes in ${\mathcal P}$). Then we again obtain a volume function defined in the neighbourhood of such a hyperplane, and can investigate the algebraicity of the analytic continuation of this function. Archimedes \cite{arch} has calculated such volumes cut from the standard (rotation-invariant) two-component hyperboloid and elliptic paraboloid in $\R^3$ by the planes orthogonal to the axes of these bodies. A trivial exercise extends this calculation to the case of non-orthogonal hyperplanes and arbitrary dimensions, and proves that 

A) the analytic continuation of the volume function on ${\mathcal P}$ defined by the two-component hyperboloid is as algebraic as the volume function defined by an ellipsoid in the space of the same dimension; 

B) for the elliptic paraboloids (i.e. the hypersurfaces having the equation 
$x_N=x_1^2+ \dots +x_{N-1}^2$
in some affine coordinates in $\R^N$) the situation is always better than for ellipsoids. Namely, for even $N$ the volume function defined by such a paraboloid is algebraic (while the one defined by an ellipsoid is not), and for odd $N$ it is polynomial while the function defined by an ellipsoid is only algebraic two-valued. 
\medskip

The fact A) has the following precise form.
\medskip

\noindent
{\bf Proposition.}
{\it Let ${\mathcal P}^{\C}$ be the space of all complex affine hyperplanes in $\C^N$; let 
$V^{\C}_W$ $($res\-pec\-ti\-vely, $V^{\C}_{W'})$ be the analytic continuation to ${\mathcal P}^{\C}$ of the volume function defined by an arbitrary ellipsoid $W \subset \R^N$ $($respectively, by an arbitrary two-component hyperboloid $W' \subset \R^N$ and any component of it$)$. Then there are an automorphism $T$ of the space ${\mathcal P}^{\C}$ and a constant $c \neq 0$ such that $V^{\C}_W \equiv c \cdot V^{\C}_{W'} \circ T$.}
\medskip

Indeed, the map $T$ is induced by the complex affine transformation of $\C^N$ moving $W$ to $W'$, and $c$ is the determinant of this transformation. \hfill $\Box$ \medskip

The ramification of these analytic continuations is defined by the Picard-Lefschetz monodromy action of the fundamental group of the space of generic affine hyperplanes in $\C^N$ on the group 
\begin{equation}
\label{hg}
H_n(\C^n, X \cup A),
\end{equation}
where $A$ is the complexification of the hypersurface $\partial W$ or $\partial W'$, and $X$ is the complexification of a generic hyperplane, see \cite{APLT}. For an explicit description of this fundamental group and this action in the case of non-singular affine conics, see \cite{Yang}. 

In the case of a paraboloid, the analogous group (\ref{hg}) is one-dimensional, the fundamental group of the space of generic affine hyperplanes in $\C^N$ (i.e. of hyperplanes neither tangent nor asymptotic to the hypersurface $A$) is isomorphic to $\Z$; its action on the group (\ref{hg}) is trivial if $N$ is odd, and is the multiplication by $-1$ if $N$ is even. This explains the fact B).

\end{document}